\documentclass[11pt,reqno,tbtags,draft]{amsart}
\title%[]
{On degenerate sums of $m$-dependent variables}

\date{5 December, 2013}
\author{Svante Janson}
\thanks{Partly supported by the Knut and Alice Wallenberg Foundation}
\address{Department of Mathematics, Uppsala University, PO Box 480,
SE-751~06 Uppsala, Sweden}
\email{svante.janson@math.uu.se}
\newcommand\urladdrx[1]{{\urladdr{\def~{{\tiny$\sim$}}#1}}}
\urladdrx{http://www2.math.uu.se/~svante/}

\subjclass[2010]{60G10; 60F05, 60C05} 
\usepackage{amssymb}
\usepackage{url}
\usepackage[square,numbers]{natbib}
\bibpunct[, ]{[}{]}{;}{n}{,}{,}
\numberwithin{equation}{section}

\renewcommand\le{\leqslant}
\renewcommand\ge{\geqslant}

\allowdisplaybreaks

\newtheorem{theorem}{Theorem}[section]

\newtheorem{corollary}[theorem]{Corollary}

\theoremstyle{definition}
\newtheorem{example}[theorem]{Example}

\newtheorem{remark}[theorem]{Remark}

\newtheorem*{ack}{Acknowledgement}

\theoremstyle{remark}

\newenvironment{romenumerate}[1][0pt]{% optional argument changes indentation
\addtolength{\leftmargini}{#1}\begin{enumerate}% gives (i), (ii) etc.
 }{\end{enumerate}}

\newcounter{oldenumi}
{\setcounter{oldenumi}{\value{enumi}}
\begin{romenumerate} \setcounter{enumi}{\value{oldenumi}}}
{\end{romenumerate}}

\newcounter{thmenumerate}

\newcounter{romxenumerate}   %less indented than standard.

\newcounter{xenumerate}   %no left indentation; thus wider lines

\newcommand{\refT}[1]{Theorem~\ref{#1}}
\newcommand{\refC}[1]{Corollary~\ref{#1}}

\newcommand{\refR}[1]{Remark~\ref{#1}}
\newcommand{\refS}[1]{Section~\ref{#1}}

\begingroup
  \divide\count255 by 60
  \multiply\count255 by -60
  \advance\count255 by \time
  \ifnum \count255 < 10 \xdef\klockan{\the\count1.0\the\count255}
  \else\xdef\klockan{\the\count1.\the\count255}\fi
\endgroup

\newcommand{\sumin}{\sum_{i=1}^n}

\newcommand\set[1]{\ensuremath{\{#1\}}}

\newcommand\bigpar[1]{\bigl(#1\bigr)}

\newcommand\biggpar[1]{\biggl(#1\biggr)}
\newcommand\lrpar[1]{\left(#1\right)}

\def\rompar(#1){\textup(#1\textup)}    % usage: \rompar(...)

\def\xexp(#1){e^{#1}}

\newcommand\ntoo{\ensuremath{{n\to\infty}}}

\newcommand\norm[1]{\|#1\|}

\newcommand\punkt{.\spacefactor=1000}    % om problem!
\newcommand\iid{i.i.d\punkt}    
\newcommand\ie{i.e\punkt}
\newcommand\eg{e.g\punkt}

\newcommand{\as}{a.s\punkt}

\newcommand{\tend}{\longrightarrow}
\newcommand\dto{\overset{\mathrm{d}}{\tend}}

\newcommand\wto{\overset{\mathrm{w}}{\tend}}
\newcommand\llto{\overset{{L^2}}{\tend}}

\newcommand\eqd{\overset{\mathrm{d}}{=}}

\newcommand\bbR{\mathbb R}

\newcommand\bbN{\mathbb N}

\newcommand\bbZ{\mathbb Z}

\newcounter{CC}
\newcounter{cc}
\newcommand\E{\operatorname{\mathbb E{}}}
\renewcommand\P{\operatorname{\mathbb P{}}}
\newcommand\Var{\operatorname{Var}}
\newcommand\Cov{\operatorname{Cov}}

\newcommand\sign{\operatorname{sign}}

\newcommand\ga{\alpha}
\newcommand\gb{\beta}

\newcommand\gs{\sigma}
\newcommand\gss{\sigma^2}

\renewcommand\phi{\xxx}  %% WARNING

\newcommand\cF{\mathcal F}

\newcommand\cT{{\mathcal T}}

\newcommand\qw{^{-1}}

\newcommand\qq{^{1/2}}

\newcommand\oi{[0,1]}

\newcommand\setoi{\set{0,1}}

\newcommand\oooo{_{-\infty}^\infty}

\newcommand\normll[1]{\norm{#1}_2}
\newcommand\bcF{\bar{\cF}}
\newcommand\xas{\text{\quad a.s.}}
\newcommand\ctn{\cT_n}
\newcommand\ct{\cT}
\newcommand\stx{\mathfrak T^*}
\newcommand\GW{Galton--Watson}
\newcommand\GWt{\GW{} tree}
\newcommand\cGWt{conditioned \GW{} tree}
\newcommand\GWp{\GW{} process}
\newcommand\tS{\tilde S}
\begin{document}

\begin{abstract} 
It is well-known that the central limit theorem holds for partial sums of a
stationary sequence $(X_i)$
of $m$-dependent random variables with finite variance; however, the limit
may be degenerate with variance 0 even if $\Var(X_i)\neq0$.
We show that this happens only in the case when $X_i-\E X_i=Y_i-Y_{i-1}$ for
an $(m-1)$-dependent stationary sequence $(Y_i)$ with finite variance
(a result implicit in earlier results), and give
a version for block factors. This yields a simple criterion that is
a sufficient condition for the limit not to degenerate. Two applications to
subtree counts in random trees are given.
\end{abstract}

\maketitle

\section{Introduction and results}\label{S:intro}

Consider a strictly stationary sequence $(X_k)_{-\infty}^\infty$ of
$m$-dependent random variables, for some $m\ge1$, 
and suppose that the variables
have finite variance, \ie, $\E X_k^2<\infty$.
(Recall that $m$-dependence means that $(X_k)_{k\le0}$ is independent of
$(X_k)_{k\ge m+1}$.)

Let $S_n:=\sumin X_i$.
A simple standard calculation using stationarity and $m$-dependence
yields, for $n\ge m$,
\begin{equation}\label{gss1}
  \begin{split}
\Var (S_n) 
&= \sum_{i,j=1}^n \Cov(X_i,X_j)
=n\Var(X_0) + 2\sum_{k=1}^m (n-k)\Cov(X_0,X_k)	
\\&
=n \gss - 2 \sum_{k=1}^m k\Cov(X_0,X_k),
  \end{split}
\raisetag{1.5\baselineskip}
\end{equation}
where
\begin{equation}\label{gss}
  \gss:= \Var(X_0) + 2\sum_{k=1}^m \Cov(X_0,X_k)	
=\Cov\biggpar{X_0,\sum_{k=-m}^m X_k}.
\end{equation}
In particular, 
\begin{equation}\label{gss3}
  \Var(S_n) = n\gss+O(1).
\end{equation}

It is obvious from \eqref{gss3} that $\gss\ge0$. If we have strict
inequality, $\gss>0$, then $\Var(S_n)$ grows linearly; moreover,
the classic central limit theorem for $m$-dependent variables 
by \citet{HoeffdingR} and \citet{Diananda},
see also \citet[Theorem 10.8]{Bradley},
shows that
\begin{equation}\label{clt}
\frac{S_n-\E S_n}{\sqrt n}=
\frac{S_n-n\E X_0}{\sqrt n} \dto N(0,\gss).
\end{equation}
In the exceptional case $\gss=0$, however, 
$\Var(S_n)$ is bounded; more precisely, \eqref{gss1} shows that $\Var(S_n)$
is constant for all $n\ge m$. In this case, \eqref{clt} still holds, with
the limit 0, but is a triviality.
(See \refC{C1} below for the limit of $S_n$ without normalization in this case.)

The purpose of the present paper is to study this exceptional case further,
and show that it really is exceptional and only occurs in very special cases.

A well-known trivial example with $\gss=0$ is obtained by taking an \iid{}
sequence 
$(Y_k)\oooo$ (with $\E Y_k^2<\infty$) and defining $X_k:=Y_k-Y_{k-1}$.
(See \eg{} \cite[\S18.1]{IL}.)
This sequence is obviously $1$-dependent and $S_n=Y_n-Y_0$ with
$\Var(S_n)=2\Var(Y_0)$, $n\ge1$, so $\Var(S_n)$ is constant and $\gss=0$. 
(This can also be seen from \eqref{gss}, using $\Var(X_0)=2\Var(Y_0)$ and
$\Cov(X_0,X_1)=-\Var(Y_0)$.) 

In fact, 
the following theorem (which is implicit in \citet[Theorem 8.6]{Bradley} but
deserves to be made more explicit)
shows that this trivial example is
the only example when $m=1$ (apart from adding a constant), and that  
a similar result holds for $m>1$. 

\begin{theorem}\label{T1}
Let $(X_k)_{-\infty}^\infty$ be a
strictly stationary sequence  of $m$-dependent variables with
finite variance and let  $\gss:=\lim_\ntoo n\qw\Var(S_n)$,
which also is given by \eqref{gss}. 
If $\gss=0$, then there exists
a strictly stationary sequence $(Y_k)_{-\infty}^\infty$ of $(m-1)$-dependent
variables with finite variance, and a constant $\mu$, such that 
\begin{equation}\label{t1}
X_k=Y_k-Y_{k-1}+\mu \xas    
\end{equation}
The random variables $Y_k$ are \as{} unique up to an additive constant.

Conversely, for any such sequence $(Y_k)\oooo$ and any $\mu$,
\eqref{t1} yields a strictly stationary $m$-dependent sequence $(X_k)\oooo$ with
$\gss=0$. 
\end{theorem}
Taking expectations in \eqref{t1} yields $\mu=\E X_k$.

\begin{remark}\label{Rweak}
 \refT{T1} holds also for weakly stationary sequences
  $(X_k)_{-\infty}^\infty$, with
  $(Y_k)_{-\infty}^\infty$ weakly stationary.
(Recall that ``weakly stationary'' just means that the means and covariances
  are translation invariant.)
\end{remark}

%\begin{remark}\label{R1}
  The existence of a (weakly)
stationary sequence $(Y_k)\oooo$ such that \eqref{t1}
  holds was shown by Leonov \cite{Leonov} under much weaker conditions
than $m$-dependence:
$(X_k)\oooo$ (weakly) stationary, $\Cov(X_0,X_n)\to0$ as \ntoo, and
$\liminf_\ntoo \Var(S_n)<\infty$. 
See also \citet{Robinson},
\citet[Theorem 18.2.2]{IL} and \citet[Theorem 8.6]{Bradley}.
The $(m-1)$-dependence of
$(Y_k)\oooo$ when $(X_k)\oooo$ is $m$-dependent 
follows from 
\cite[Theorem 8.6(B)(e)]{Bradley}, but does not seem to have been stated
explicitly earlier.

For completeness, we give a direct proof of \refT{T1} in \refS{Spf1}.
(The same proof applies to the weakly stationary version, see \refR{Rweak}.)

\begin{remark}
  More generally, a theorem by \citet[Lemma 11.7]{Schmidt},
in the version given by \citet[Theorem 19.9]{Bradley},
implies that even without the assumption of finite variance, 
if $(X_k)\oooo$ is a strictly stationary and $m$-dependent
sequence such that the family of partial sums $S_n$ are tight,
then the conclusion
\eqref{t1} (with $\mu=0$) holds for some strictly stationary
$(m-1)$-dependent sequence $(Y_k)\oooo$. 
\end{remark}

Note that \eqref{t1} implies
\begin{equation}\label{sn}
S_{n}-\E S_n = S_n -n\mu =Y_n-Y_{0} \xas,
\end{equation}
where $Y_n\eqd Y_0$ and  $Y_n$ and $Y_0$ are independent when $n\ge m$.
%In particular, $\Var(S_n)=2\Var(Y_0)$ when $n\ge m$.
An immediate consequence of \refT{T1} is thus that in the exceptional case
$\gss=0$, the centered partial sums $S_n-\E S_n$ converge in distribution
without normalization. Of course, the limit is in general not normal, so
there is no central limit theorem in this case. (For example, $X_n$ may be
integer valued, and then so is $S_n$.) We state this in detail.
(See \refS{Spf1} for proofs of this and other results.)

\begin{corollary}\label{C1}
Let $(X_k)_{-\infty}^\infty$ be a
strictly stationary sequence  of $m$-dependent variables with
finite variance, and let $\gss$
be given by \eqref{gss}.
If $\gss=0$, then $S_n-\E S_n$ has the same distribution for all  
$n\ge m$; more precisely, if $Y_k$ is as in \eqref{t1} and $Y_0'$ is an
independent copy of $Y_0$, then $S_n-\E S_n\eqd Y_0-Y'_0$.

Hence, assuming that $\Var(X_0)>0$, $(S_n-\E S_n)/\Var(S_n)\qq$ converges in
distribution as \ntoo{} also in the case $\gss=0$, but then the limit is 
normal only if each $Y_k$ is normal.
\end{corollary}

\begin{remark}\label{RN}
  If $\gss=0$ and $m=1$, then \eqref{t1}
holds with independent $Y_k$. Hence, by a
theorem by Cram{\'e}r, 
see \eg{} \cite[Theorem XV.8.1]{FellerII},
each $Y_k$ is normal if and
only if $X_k$ is normal (and then $\set{X_k,Y_k:k\in\bbZ}$
are jointly normal). For
$m>1$ this does not hold. For example,  
if $U_k\sim U(0,1)$ and $\xi_k\sim N(0,1)$, $k\in\bbZ$, 
all  independent, then
$Y_k:=\sign(U_k-U_{k+1})|\xi_k|$ is a sequence of 1-dependent normal
variables that are not jointly normal, and  
the 2-dependent random variables $X_k:=Y_k-Y_{k-1}$ 
are not
normal although by \eqref{sn}  $S_n\sim N(0,2)$ is for $n\ge2$.
(A simple calculation yields $\E X_k^2=2+\frac{4}{3\pi}$ and 
$\E X_k^4=12+\frac{32}{3\pi} \neq 3 (\E X_k^2)^2$.)
\end{remark}

Stationary
$m$-dependent sequences usually appear as \emph{block factors}.
We say that $(X_k)$ is an $\ell$-block factor if there is an \iid{} sequence
$(\xi_k)\oooo$ and a (measurable) function $f:\bbR^\ell\to\bbR$ such that
$X_k=f(\xi_k,\dots,\xi_{k+\ell-1})$.
Note that every such sequence $(X_k)$ is strictly stationary and
$(\ell-1)$-dependent.
(However, there are $m$-dependent sequences that are not block factors
\cite{Aaronson}, \cite{Burton}.)

For block factors, \refT{T1} takes the following form.

\begin{theorem}\label{T2}
Let $X_k=f(\xi_k,\dots,\xi_{k+\ell-1})$ be an $\ell$-block factor
for some $\ell\ge1$,
where $(\xi_k)_{-\infty}^\infty$ is an \iid{} sequence.
Suppose that $X_k$ has
finite variance and let  $\gss:=\lim_\ntoo n\qw\Var(S_n)$.
If $\gss=0$, then there exists a function $g:\bbR^{\ell-1}\to\bbR$ and a
constant $\mu$ 
such that the $(\ell-1)$-block factor
$Y_k:=g(\xi_{k+1},\dots,\xi_{k+\ell-1})$ has finite variance and
\begin{equation}\label{t2}
X_k=Y_k-Y_{k-1}+\mu \xas    
\end{equation}
The function $g$ is \as{} unique up to an additive constant.
%Conversely, for any such sequence $(Y_k)\oooo$ and any $\mu$,
%\eqref{t1} yields a strictly stationary $m$-dependent sequence $(X_k)\oooo$
%with $\gss=0$. 
\end{theorem}
The converse is obvious in this theorem too.
%Moreover, the theorem leads to a simple criteris fo

\begin{corollary}
  \label{C2}
Let $X_k=f(\xi_k,\dots,\xi_{k+\ell-1})$ be an $\ell$-block factor
with finite variance,
where $(\xi_k)_{-\infty}^\infty$ is an \iid{} sequence.
If $\gss=0$, then there exists a function $g:\bbR^{\ell-1}\to\bbR$ 
such that for every $n\ge1$,
\begin{equation}\label{c2}
S_n-\E S_n=g(\xi_{n+1},\dots,\xi_{n+\ell-1})
-g(\xi_{1},\dots,\xi_{\ell-1}) \xas
\end{equation}
\end{corollary}

\begin{remark}\label{RC2}
The contrapositive form of \refC{C2} yields a simple criterion:
If we can find, for some $n\ge \ell$,
a set of values of $\xi_1,\dots,\xi_{\ell-1}$
and $\xi_{n+1},\dots,\xi_{n+\ell-1}$ of positive probability such that
$S_n$ is not an \as{} constant function of $\xi_{\ell},\dots,\xi_{n}$,
then \eqref{c2} cannot hold and thus $\gss>0$.
\end{remark}

\refC{C2} and its reformulation in \refR{RC2}
are useful in applications, 
to show that the asymptotic variance $\gss>0$.
We give two such applications in \refS{Sex},
taken from \citet{HolmgrenJ} and \citet{SJ285};
these applications were the motivation for the present study.

\begin{remark}
  The central limit theorem for $m$-dependent variables has been generalized
  to much more general mixing sequences under various conditions,
see \eg{} \cite{IL} and \cite{Bradley}.
For example, if $(X_k)\oooo$ is strictly stationary 
with finite variances
and $\rho$-mixing, 
then
either 
\begin{romenumerate}
\item 
$\Var(S_n)=nh(n)$ for some slowly varying function $h(n)$, or
\item $\Var(S_n)$ is bounded, and converges to some finite limit.
\end{romenumerate}
Moreover, in case (i), a central limit theorem holds under 
further conditions
\cite{Ibragimov}, 
\cite[Theorems 11.2 and 11.4]{Bradley}
(but not in general 
\cite{Bradley80},
\cite[Chapter 34]{Bradley}).

In case (ii), there is by the result by \citet{Leonov} mentioned above 
a representation as in
 \eqref{t1}; however,
we do not know any useful consequences similar to \refC{C2} and \refR{RC2}
in this generality and 
we leave it as  an open problem to find generalizations of the results above
that can be used to show $\gss>0$.
A typical example of case (ii) is $X_k=\xi_k-\sum_{j=1}^\infty 2^{-j}\xi_{k+j}$
with $(\xi_k)\oooo$ \iid{} $N(0,1)$,
where we have the representation \eqref{t1} with
$Y_k=-\sum_{j=0}^\infty 2^{-j}\xi_{k+1+j}$.
\end{remark}

\begin{ack}
  I thank Richard Bradley for several valuable comments.
\end{ack}

\section{Proofs}\label{Spf1}

\begin{proof}[Proof of \refT{T1}]
As said in the introduction, \refT{T1} follows from \cite[Theorem 8.6]{Bradley},
but we give also a direct proof for completeness. (The proof is similar, 
but  simpler in this special case.)

It is obvious that if $(Y_k)\oooo$ is strictly stationary and
$(m-1)$-dependent,
then $(X_k)\oooo$ defined by \eqref{t1} is strictly stationary and
$m$-dependent. 
Furthermore, \eqref{t1} implies \eqref{sn} and thus 
$\Var(S_n)=\Var(Y_n)+\Var(Y_0)=2\Var(Y_0)$ when $n\ge m$; hence
$\gss=0$ by \eqref{gss3}.

To prove the converse we may assume $\E X_k=0$.
  Define 
$S_{k,n}:= \sum_{i=k}^n X_i$, for $-\infty <k\le n<\infty$.
The assumption $\gss=0$ implies by \eqref{gss3} and stationarity that
$\E S_{k,n}^2=\Var(S_{k,n})$ is bounded. 
(In fact, by \eqref{gss1} it is constant for all $(k,n)$
with $n-k\ge m-1$.) 

We claim first that for every $k$, the sequence $S_{k,n}$ converges weakly in
$L^2$ as \ntoo, and thus there exists a random variable $Z_k\in L^2$ such
that
\begin{equation}\label{wto}
  S_{k,n} \wto Z_k \qquad \text{as }\ntoo.
\end{equation}
In fact, since the sequence $(S_{k,n})_{n\ge k}$ is bounded in $L^2$ and
the unit ball of $L^2$ is weakly compact, it suffices to show that 
$\E(W  S_{k,n})$ converges as \ntoo{} for every fixed $W\in L^2$; moreover,
it suffices to verify this for a dense set of $W$. 
We consider two special cases:
\begin{romenumerate}
\item 
If $\E (W X_j)=0$ for
all $j$, then $\E(W S_{k,n})=0$ for all $n$, and the convergence is trivial.
\item 
 If $W=X_j$ for some $j$, then 
$\E(W S_{k,n})$ is constant for all $n\ge \max(j+m,k)$, by $m$-dependence,
and again the convergence is trivial. 
\end{romenumerate}
Hence $\E(W  S_{k,n})$ converges also when $W$ is a linear combination of
variables of the type (i) or (ii). But the set of such linear combinations
is dense in $L^2$, which proves \eqref{wto}.

Similarly (or by reflecting the indices and replacing $X_k$ by $X_{-k}$),
for every $k\in\bbZ$ there exists a random variable $Y_k\in L^2$ such that
\begin{equation}\label{wtoy}
  S_{-n,k}\wto Y_k \qquad \text{as }\ntoo.
\end{equation}
Since $S_{-n,k} - S_{-n,k-1} = X_k$ for $-n<k$, it follows that
$Y_k-Y_{k-1}=X_k$, so \eqref{t1} holds (with $\mu=\E X_0=0$).
Furthermore, $(Y_k)\oooo$ is stationary by
\eqref{wtoy} and the stationarity of $(X_k)\oooo$.
It remains to show that $(Y_k)\oooo$ is $(m-1)$-dependent.

We note first that for any $k$, as \ntoo,
by \eqref{wtoy} and \eqref{wto},
\begin{equation}\label{ss1}
  S_{-n,k}+S_{k+1,n} \wto Y_k+Z_{k+1}.
\end{equation}
On the other hand,
$S_{-n,k}+S_{k+1,n}= S_{-n,n}$ (when $n>|k|$), and 
thus for every $j\in \bbZ$ and every $n
> \max(|k|,m+|j|)$,
using $m$-dependence and \eqref{gss}, 
\begin{multline}\label{ss2}
\E\bigpar{X_j( S_{-n,k}+S_{k+1,n})}
=\E(X_j S_{-n,n}) =\Cov(X_j,S_{-n,n})
\\
\hskip-4em
=\sum_{i=-n}^n\Cov(X_j,X_i)
=\sum_{i=j-m}^{j+m}\Cov(X_j,X_i)	
=\gss
=0.	  
\end{multline}
Combining \eqref{ss1} and \eqref{ss2} we see that $\E(X_j (Y_k+Z_{k+1}))=0$
for every $j$. Summing over $j$ we find
$\E(S_{\ell,n} (Y_k+Z_{k+1}))=0$ for all $\ell$ and $n$, and thus
by \eqref{ss1} again,
$\E(Y_k+Z_{k+1})^2=0$. Hence $Y_k+Z_{k+1}=0$ a.s., \ie{}
\begin{equation}\label{yz}
  Y_k=-Z_{k+1} \xas
\end{equation}

For $-\infty\le k\le n\le\infty$,
let  $\cF_{k,n}$ denote the $\gs$-field generated by $\set{X_i}_{i=k}^n$.
Write $W\in\cF_{k,n}$ if the random variable $W$ is $\cF_{k,n}$-measurable.
Then $S_{-n,k}\in \cF_{-n,k}\subseteq\cF_{-\infty,k}$, and thus \eqref{wtoy}
shows that 
\begin{equation}\label{yk-}
Y_{k}\in \cF_{-\infty,k}.  
\end{equation}
Similarly, $Z_{k}\in \cF_{k,\infty}$.  By \eqref{yz}, this yields also
\begin{equation}\label{yk+}
Y_{k}\in \cF_{k+1,\infty}.  
\end{equation}

Since $(X_k)\oooo$ is $m$-dependent, the $\gs$-fields
$\cF_{-\infty,k}$ and $\cF_{k+m+1,\infty}$ are independent. Hence,
\eqref{yk-}--\eqref{yk+} show that
\set{Y_j:j\le k} and $\set{Y_j:j\ge k+m}$ are independent, for every $k$,
which is the desired ($m-1$)-dependence.

Finally, we consider uniqueness of $Y_k$. It is obvious that we may replace
$Y_k$ by $Y_k+C$ for any constant $C$.
For the converse, we may assume $\E X_k=0$ so $\mu=0$.
If \eqref{t1} holds, then
\begin{equation}\label{skn}
S_{k,n}=Y_n-Y_{k-1},
\end{equation}
and it follows by \eqref{gss3} applied to $(Y_n)\oooo$ that
\begin{equation}\label{cesaro}
  \Var 
\biggpar{\frac1{n} \sum_{j=k+1}^{k+n}S_{k,j}+Y_{k-1} } =O\bigpar{n\qw}
\end{equation}
and thus $Y_{k-1}-\E Y_{k-1}$ is the limit in $L^2$ of the means
$-\frac1{n} \sum_{j=k+1}^{k+n}S_{k,j}$, and thus \as{} determined by
$(X_j)\oooo$. 
\end{proof}

\begin{remark}
  We use weak convergence in $L^2$ in
\eqref{wto} and \eqref{wtoy}, following Leonov
  \cite{Leonov} who uses weak convergence of a subsequence in a much more
  general situation. (It is easy to modify the proof by Leonov
  \cite{Leonov} to show weak convergence of the full sequence under the
  conditions there too. We have above used a simpler version for the
  $m$-dependent case.)
Strong (norm) convergence does not hold:
\eqref{skn} shows that
$\normll{S_{k,n}-Z_k}   = \normll{S_{k,n}+Y_{k-1}}  
= \normll{Y_{n}}$
which is constant and does not tend to 0 
(except in the trivial case $Y_n=0$ when $X_k=0$ \as{}).
However, assuming $\E X_k=0$, \eqref{cesaro} shows that the Ces{\`a}ro means 
$T_{k,n}:=(n+1)\qw \sum_{j=k}^{k+n} S_{k,j}$ 
converge to $Z_k=-Y_{k-1}$ in $L^2$, \ie,
$\normll{T_{k,n}-Z_k}\to0$, %:=\bigpar{\E(T_{n,k}-Z_k)^2}\qq\to0$, 
and similarly
$(n+1)\qw \sum_{j=k-n}^{k} S_{k,j} \llto Y_k$,
see the proof of
\cite[Theorem 8.6]{Bradley}. 
(This can be used to give an alternative proof
of \refT{T1}, using strong Ces{\`a}ro convergence instead of weak convergence
and completing the proof as above.)
Furthermore, the strong law of large numbers for stationary
$m$-dependent sequences implies $T_{k,n}\to Z_k$ a.s., while \eqref{skn}
shows that $S_{k,n}$ does not converge \as{} (except when $Z_i=0$).
\end{remark}

\begin{proof}[Proof of \refC{C1}]
  By \refT{T1}, \eqref{t1} holds and thus 
\eqref{sn} holds, which shows $S_n-\E S_n\eqd Y_0-Y_0'$ when $n\ge m$.
In particular, for $n\ge m$, 
$\Var(S_n)=2\Var(Y_0)$ and hence $\Var(S_n)=0$ only if $Y_0$ is
  degenerate (\as{} constant), and then each $X_k$ is degenerate.
Finally, by the theorem by Cram{\'e}r mentioned in \refR{RN},
%see \eg{} \cite[Theorem XV.8.1]{FellerII},
$Y_0-Y'_0$ is normal if and only if
$Y_0$ has a normal distribution.
\end{proof}

\begin{proof}[Proof of \refT{T2}]
Let  $Y_k$ and $Z_k$ be as in the proof of \refT{T1}.

For $-\infty\le k\le n\le\infty$,
let  $\bcF_{k,n}$ denote the $\gs$-field generated by $\set{\xi_i}_{i=k}^n$
and all sets of probability 0. (The latter technicality is because $Y_k$ and
$Z_k$ are defined only a.s.)
Then $X_k\in \bcF_{k,k+\ell-1}$
so $S_{k,n}\in \bcF_{k,n+\ell-1}$ and thus
$Y_{k}\in \bcF_{-\infty,k+\ell-1}$ and 
$Z_{k}\in \bcF_{k,\infty}$. Since $Y_k=-Z_{k+1}$ by \eqref{yz}, thus
\begin{equation}\label{yg}
  Y_k \in  \bcF_{-\infty,k+\ell-1}\cap \bcF_{k+1,\infty}
= \bcF_{k+1,k+\ell-1},
\end{equation}
where the latter equality follows (\eg{} by considering conditional
expectations) because the variables $\xi_i$ are independent.

Hence, $Y_k=g(\xi_{k+1},\dots,\xi_{k+\ell-1})$ for some function $g$
(independent of $k$ because of stationarity).
The result now follows from \refT{T1}.
\end{proof}

\begin{proof}[Proof of \refT{C2}]
An immediate consequence of \refT{T2} and \eqref{sn}.
\end{proof}

\section{Applications}\label{Sex}
We sketch two applications of the results above;
more details and background are given in \citet{HolmgrenJ} and
\citet{SJ285}.
In both applications we
consider a random rooted tree $\ctn$ with $n$ nodes (with different
distributions in the two cases) and let for a fixed rooted tree $T$,
$n_T(\ctn)$ be 
the number of nodes $v\in\ctn$ such that the fringe subtree consisting of
$v$ and all its descendants is isomorphic to $T$. 
(We consider only trees $T$ 
in the family $\stx$ of trees
that can appear as fringe subtrees in $\ctn$ for 
some $n$; otherwise $n_T(\ctn)$ is identically 0 for all $n$.)
In
the cases studied here, these numbers are asymptotically normal for fixed
$T$ as \ntoo:
\begin{equation}\label{ntt}
  \frac{n_T(\ctn)-n\mu_T}{\sqrt n} \dto \zeta_T
\end{equation}
where $\zeta_T\sim N\bigpar{0,\gss_T}$,
for some $\mu_T>0$ and $\gss_T\ge0$; moreover,
this holds jointly for all $T$ with the limit variables $\zeta_T$ jointly
normal, with convergence of variances and covariances.
We use the results above to show that the limit distribution is not
degenerate: $\gs^2_T>0$ for each $T\in\stx$, and moreover, the covariance
matrix of 
$\zeta_{T_1},\dots,\zeta_{T_N}$ is positive definite, for 
any finite number of trees $T_1,\dots,T_N\in\stx$.
Equivalently, 
if
\begin{equation}\label{fan}
  F(\ctn) = \sum_{i=1}^N a_jn_{T_j}(\ctn)
\end{equation}
for some distinct trees 
$T_1,\dots,T_N\in\stx$ and real numbers $a_1,\dots,a_N$, not all zero,
then
\begin{equation}\label{flim}
 \lim_\ntoo\frac{ \Var F(\ctn)} n 
=\Var\biggpar{ \sum_{j=1}^N a_j \zeta_{T_j}} >0.
\end{equation}

\begin{example}[Binary search trees, \cite{HolmgrenJ}]
A \emph{binary search tree} is a binary tree with a \emph{key} stored at
each node. It is constructed from a sequence of (distinct) keys 
by putting the
first key, say $x_1$, in the root and sending all subsequent keys less than
$x_1$ to the left subtree and the keys greater than $x_1$ to the right
subtree, constructing the subtrees recursively in the same way.

We may assume that the keys are $1,\dots,n$; then, a binary search tree 
is a binary tree with the nodes labelled $1,\dots,n$ (where $n$ is the size
of the tree).
Let $\ctn$ be a uniformly random binary search tree with $n$ nodes; this can
be constructed by taking the keys $1,\dots,n$ in (uniformly) random order, 

We use a modification of this construction by
Devroye \cite{Devroye1,Devroye2}:
Let $U_1,\dots,U_n$ be \iid{} random variables with $U_i\sim U(0,1)$,
order the indices $1,\dots,n$ so that the variables $U_i$ are in increasing
order and construct the binary search tree $\ctn$ as above using this sequence.
(Thus, for example, the root is labelled by the index $i$ such that $U_i$ is
the smallest of $U_1,\dots,U_n$.)
It is not difficult to see that then
the fringe subtrees of $\ctn$ are the trees 
defined in the same way by the subsequences $U_i,\dots,U_j$ (with $1\le i\le
j\le n$) such that $U_{i-1}$ and $U_{j+1}$ both are smaller than all of
$U_i,\dots,U_j$; we here define $U_0=U_{n+1}=0$.

Hence, if $T\in\stx$, where now $\stx$ is the family of all binary trees,
and $T$ has $|T|=k$ nodes, then  
\begin{equation}\label{bbol}
  n_{T}(\ctn) = \sum_{i=0}^{n-k} f_T(U_{i},\dots,U_{i+k+1})
\end{equation}
for some indicator function $f_T(x_1,\dots,x_{k+2})$ on $\oi^{k+2}$
(depending only on the order relations between $x_1,\dots,x_{k+2}$).
For convenience we ignore the boundary terms in \eqref{bbol}, which are
asymptotically negligible; we let
$(U_i)\oooo$ be \iid{} with $U_i\sim U(0,1)$ and then
\begin{equation}\label{bbol-}
  n_{T}(\ctn) = \sum_{i=1}^{n-k-1} f_T(U_{i},\dots,U_{i+k+1}) +O(1),
\end{equation}
where the sum is 
a sum of $m$-dependent variables of the type studied in this paper.
Given a function $F$ as in \eqref{fan}, we 
let $\ell:=\max_j|T_j|+2$ and
define $f(x_1,\dots,x_\ell):=\sum_j a_j f_{T_j}(x_1,\dots,x_{|T_j|+2})$;
then \eqref{bbol-} implies
\begin{equation}\label{gul}
  F(\ctn)= \sum_{i=1}^{n-\ell} f(U_{i},\dots,U_{i+\ell-1}) +O(1)
=S_{n-\ell}+O(1),
\end{equation}
where $S_n=\sum_{i=1}^n X_i$ with $X_i=f(U_i,\dots,U_{i+\ell-1})$ an
$\ell$-block factor as in \refT{T2}.
Hence the central limit theorem for $m$-dependent variables 
\cite{HoeffdingR}, \cite{Diananda} yields asymptotic normality of $F(\ctn)$,
\ie, \eqref{ntt} with joint convergence for several $T\in\stx$ and
convergence of first and second moments; this is the method 
by \citet{Devroye1}.
We can now also show that \eqref{flim} holds. 

We may suppose that $a_1,\dots,a_N$ all are non-zero, and that
$T_1,\dots,T_N$ are ordered with $|T_1|\le |T_2|\le\dots$, so no $T_j$ is a
proper subtree of $T_1$. Let $n>3\ell$, and consider the event that
$U_1<U_2<\dots<U_n$; this generates a tree $\ctn=T'$ that is a path to the
right from the root. By permuting $U_\ell,\dots,U_{\ell+k}$, where
$k=|T_1|$, leaving all other $U_i$ unchanged, we may instead
generate a tree $T''$ that 
is a path to the right of length $n-k$, with a copy of $T_1$
attached to the $\ell$:th vertex. Then $n_{T_1}(T'')=n_{T_1}(T')+1$,
but $n_{T_j}(T'')=n_{T_j}(T')$ for $2\le j\le N$, since 
except for the new copy of $T_1$ in $T''$, 
the subtrees
that appear or disappear when we change $T'$ to $T''$ 
are either too
small or too large to be a $T_j$. Hence, by \eqref{fan}, $F(T')\neq F(T'')$,
and this holds also if we ignore the boundary trees and consider $S_n$ as in
\eqref{gul}, and it follows by \refC{C2}, see \refR{RC2}, that \eqref{flim}
holds. 
(The proof just given was our first proof that $\gss>0$ in this case.
The proof given in \citet{HolmgrenJ} is actually slightly different
 and does not use the results in the present paper; 
it uses instead a shortcut based on a special symmetry property.)
\end{example}

\begin{example}[Conditioned Galton--Watson trees, \cite{SJ285}]
A \GWt{} $\ct$ is the tree version of a \GWp{}. It is defined by a 
non-negative integer-valued random variable $\xi$ which describes the number
of children of each node. We assume that $\E \xi=1$ (a \emph{critical} \GWp) 
and $\E\xi^2<\infty$.
The \cGWt{} $\ctn$ is the random tree $\ct$ conditioned to have exactly $n$
nodes. It is well-known that several standard types of random trees can be
defined in this way, with suitable $\xi$, see \eg{} \cite{SJ264}.
We assume for simplicity that $\P(\xi=k)>0$ for every $k\ge0$, and let
$\stx$ be the family of all ordered rooted trees.
(The general case is studied in \cite{SJ285} with a minor variation of the
argument below. The result is the same as long as $\xi$ attains at least two
positive integers with positive probability, 
except that $\stx$ only consists of trees where all outdegrees may be
attained by $\xi$,
but in the case when
$\xi\in\set{0,r}$ for some integer $r$, we have to exclude the case
$T=\bullet$, the tree of size 1, because $n_\bullet(\ctn)$ is deterministic.)

Let $\xi_1,\xi_2,\dots$ be an \iid{} sequence of copies of $\xi$, and let
$Z_n:=\sum_{i=1}^n \xi_i$. The degree sequence of the nodes in $\ctn$, taken
in depth-first order, is $(\xi_1,\dots,\xi_n)$ conditioned on this being the
degree sequence of a tree; up to a cyclic shift this is the same as
conditioning on $Z_n=n-1$ and it follows that
\begin{equation}\label{nd}
  n_{T}(\ctn)
\eqd \lrpar{\sum_{i=1}^{n} f_T(\xi_i,\dots,\xi_{i+k-1\bmod n})
\Bigm| Z_n=n-1},
\end{equation}
for a suitable indicator function $f_T:\bbN^k\to\setoi$, where $k=|T|$.
Given $F$ as in \eqref{fan}, we let $\ell:=\max_j|T_j|$ and
$f(x_1,\dots,x_\ell):=\sum_j a_j f_{T_j}(x_1,\dots,x_{|T_j|})$;
then, again ignoring some boundary terms,
\begin{equation}\label{Fd}
  \begin{split}
  F(\ctn)
&\eqd \lrpar{\sum_{i=1}^{n-\ell} f(\xi_i,\dots,\xi_{i+\ell-1})
\Bigm| Z_n=n-1} +O(1)
\\&
=\bigpar{S_{n}\mid Z_n=n-1} + O(1).	
  \end{split}
\end{equation}
In this case, we thus have a conditioned version of the sum
$S_n$, and asymptotic normality follows by a method by \citet{LeCam} and
\citet{Holst}, see \cite{SJ285}. The proof shows that the asymptotic
variance $\gss$ is given by
\begin{equation}
  \gss = \lim_{\ntoo}\frac{1}{n}\Var\bigpar{S_n-\ga Z_n},
\end{equation}
where the constant $\ga$ is chosen such that 
$\Cov\bigpar{S_n-\ga Z_n,Z_n}/n\to0$.
Let  $\tS_n:=S_n-\ga Z_n$. Then $\tS_n-\E\tS_n = \sumin X_i$, where
\begin{equation}\label{Xi}
  X_i := f(\xi_i,\dots,\xi_{i+\ell-1}) -\ga \xi_i + \gb,
\end{equation}
with $\gb$ chosen such that $\E X_i=0$.
If $\gss=0$, we may apply \refC{C2}. 
Take first $\xi_i=j$ for all $i\le n+\ell-1$, for some $j>0$.
Then $(\xi_i,\dots,\xi_{i+k-1})$ is never the degree sequence of a tree, so 
$f_T(\xi_i,\dots,\xi_{i+k-1})=0$  and 
$f(\xi_i,\dots,\xi_{i+\ell-1})=0$; hence \eqref{Xi} reduces to 
$X_i=-\ga j+\gb$, and \eqref{c2} yields $n(-\ga j+\gb)=0$.  
Hence $-\ga j + \gb=0$ for every $j>0$, and thus $\ga=\gb=0$.
We may again assume that $|T_1|\le|T_2|\le\dots\le |T_N|$ and $a_1\neq0$.
Let $n>2\ell$ and assume that $(\xi_{\ell+1},\dots,\xi_{\ell+|T_1|})$
equals the degree sequence of $T_1$, while all other  $\xi_i=2$, say, for
$i\le n+\ell-1$. The only substrings of $\xi_1,\dots,\xi_{n+\ell-1}$ that
are degree sequences of trees are 
$(\xi_{\ell+1},\dots,\xi_{\ell+|T_1|})$ and some of its substrings,
corresponding to $T_1$ and its subtrees. It follows that 
$\tS_n-\E \tS_n = a_1\neq0$, which contradicts \eqref{c2}. This
contradiction proves $\gss>0$, \ie, \eqref{flim}.
\end{example}

\newcommand\AAP{\emph{Adv. Appl. Probab.} }
\newcommand\JAP{\emph{J. Appl. Probab.} }
\newcommand\JAMS{\emph{J. \AMS} }
\newcommand\MAMS{\emph{Memoirs \AMS} }
\newcommand\PAMS{\emph{Proc. \AMS} }
\newcommand\TAMS{\emph{Trans. \AMS} }
\newcommand\AnnMS{\emph{Ann. Math. Statist.} }
\newcommand\AnnPr{\emph{Ann. Probab.} }
\newcommand\CPC{\emph{Combin. Probab. Comput.} }
\newcommand\JMAA{\emph{J. Math. Anal. Appl.} }
\newcommand\RSA{\emph{Random Struct. Alg.} }
\newcommand\ZW{\emph{Z. Wahrsch. Verw. Gebiete} }
\newcommand\DMTCS{\jour{Discr. Math. Theor. Comput. Sci.} }

\newcommand\AMS{Amer. Math. Soc.}
\newcommand\Springer{Springer-Verlag}
\newcommand\Wiley{Wiley}

\newcommand\vol{\textbf}
\newcommand\jour{\emph}
\newcommand\book{\emph}
\newcommand\inbook{\emph}
\def\no#1#2,{\unskip#2, no. #1,} %(typeset after year) 
\newcommand\toappear{\unskip, to appear}

\newcommand\urlsvante{\url{http://www.math.uu.se/~svante/papers/}}
\newcommand\arxiv[1]{\url{arXiv:#1.}}
\newcommand\arXiv{\arxiv}

\def\nobibitem#1\par{}

\end{document}